\documentclass{article}
\usepackage{amsmath,amsfonts,amsthm,amssymb,amscd,bm,bbm}
\usepackage[hidelinks]{hyperref}
\usepackage[usenames,dvipsnames]{color}
\usepackage{url}
\hypersetup{colorlinks, citecolor=ForestGreen, linkcolor=MidnightBlue, urlcolor=Blue}

\newcommand{\be}{\begin{equation}}
\newcommand{\ee}{\end{equation}}
\newcommand{\bp}{\begin{proof}}
\newcommand{\ep}{\end{proof}}
\newcommand{\bi}{\begin{itemize}}
\newcommand{\ei}{\end{itemize}}
\newcommand{\om}{\omega}

\newcommand{\iI}{\mathfrak{i}}

\newcommand{\bB}{\mathfrak{b}}
\newcommand{\uu}{\mathfrak{u}}

\newcommand{\aaa}{{\cal A}}
\newcommand{\BBB}{\mathfrak B}

\newcommand{\ees}{{\cal E}}
\newcommand{\fff}{{\cal F}}

\newcommand{\mmm}{{\cal M}}

\newcommand{\ooo}{{\cal O}}
\newcommand{\ppp}{{\cal P}}

\newcommand{\R}{{\cal R}}

\newcommand{\vvv}{{\cal V}}

\newcommand{\yyy}{{\cal Y}}

\newcommand{\e}{\mathbb{E}}

\newcommand{\nn}{\mathbb{N}}
\newcommand{\pp}{\mathbb{P}}
\newcommand{\qq}{\mathbb{Q}}
\newcommand{\rr}{\mathbb{R}}

\newcommand{\q}{\quad}

\newcommand{\f}{\frac}
\newcommand{\lm}{\lambda}
\newcommand{\p}{\partial}
\newcommand{\ph}{\varphi}

\newcommand{\De}{\delta}
\newcommand{\de}{\Delta}
\newcommand{\g}{\nabla}
\newcommand{\dt}{\dot}
\newcommand{\diver}{\mathop{\rm div}\nolimits}

\newcommand{\es}{\varepsilon}

\newcommand{\al}{\alpha}

\newcommand{\nt}{\noindent}

\newcommand{\ch} {\mathbbm{1}}

\newcommand{\iin} {\infty}

\newcommand{\mpp}{\emph}
\newcommand{\ef}{\eqref}
\newcommand{\dd}{\,{\textup d}}
\newcommand{\Dd}{{\textup d}}
\theoremstyle{plain}
\newtheorem{theorem}{Theorem}[section]

\newtheorem{lemma}[theorem]{Lemma}
\newtheorem{proposition}[theorem]{Proposition}

\theoremstyle{definition}
\newtheorem{definition}[theorem]{Definition}

\theoremstyle{remark}

\numberwithin{equation}{section}
\newtheorem*{definition*}{Definition}
\newtheorem*{problem*}{Problem}

\newtheorem*{remark*}{Remark}
\newtheorem*{note*}{Note}
\begin{document}
\author{Davit Martirosyan\footnote{Department of Mathematics, University of Cergy-Pontoise, CNRS UMR 8088, 2 avenue
Adolphe Chauvin, 95300 Cergy-Pontoise, France;e-mail: \href{mailto:Davit.Martirosyan@u-cergy.fr}{Davit.Martirosyan@u-cergy.fr}}}
\title{Large deviations for invariant measures of white\,-\,forced 2D Navier-Stokes equation} 
\date{}
\maketitle

\begin{center}
{\it To the memory of William Unterwald (adp)}
\end{center}
\begin{abstract}
The paper is devoted to studying the asymptotics of the family $(\mu^\es)_{\es>0}$ of stationary measures of the Markov process generated by the flow of equation
$$
\dt u-\de u+(u,\g)u+\g p=h(x)+\sqrt{\es}\,\vartheta(t,x), \q \diver u=0, \q u|_{\p D}=0
$$
in a bounded domain $D\subset \rr^2$, where $h\in H^1_0(D)$ and $\vartheta$ is a spatially regular white noise. By using the large deviations techniques, we prove that the family $(\mu^\es)$ is exponentially tight in $H^{1-\gamma}(D)$ for any $\gamma>0$ and vanishes exponentially outside any neighborhood of the set $\ooo$ of $\om$-limit points of the deterministic equation. In particular, any of its weak limits is concentrated on the closure $\bar\ooo$. A key ingredient of the proof is a new formula that allows to recover the stationary measure $\mu$ of a Markov process with good mixing properties, knowing only some local information about $\mu$. In the case of trivial limiting dynamics, our result implies that the family $(\mu^\es)$ obeys the large deviations principle. 

 \smallskip
\noindent
{\bf AMS subject classifications:}  35Q30, 76D05, 60F10, 60H15

\smallskip
\noindent
{\bf Keywords:} large deviations principle, stochastic partial differential equations, invariant measures, white noise
\end{abstract}

\tableofcontents

\setcounter{section}{-1}

\section{Introduction}\label{1.77}
Let us consider the Navier-Stokes equation
 \be\label{0.61}
\dt u-\de u+(u,\g )u+\g p=f(t,x), \q \diver u=0, \q u|_{\p D}=0
\ee
in a bounded domain $D\subset \rr^2$.  As is well-known, \ef{0.61} gives rise to an evolution equation if we eliminate the pressure term using the Leray projection $\Pi: L^2(D)\to H$, where $H$ is the space of divergence-free vector fields of $L^2(D)$ with vanishing normal component. The corresponding equation reads
$$
\dt u+Lu+B(u,u)=\Pi f(t,x),
$$
where $L$ is the Stokes operator and $B(u, v)$ stands for the bilinear form $\Pi(u, \g)v$.

\medskip
In this paper we study the asymptotics of the family $(\mu^\es)_{\es>0}$ of stationary measures of the Markov process generated by the flow of equation
\be\label{0.8}
\dt u+ Lu+B(u,u)=h(x)+\sqrt{\es}\,\vartheta(t,x),
\ee
where $h$ is a function in $H^1_0(D)\cap H$ and $\vartheta$ is a colored white noise of the form
\be
\vartheta(t,x)=\sum_{m=1}^\iin b_m\dt \beta_m(t)e_m(x).
\ee
Here $\{\beta_m(t)\}$ is a sequence of standard Brownian motions, $\{e_m\}$ is an orthogonal basis in $H$ composed of the eigenfunctions of the Stokes operator $L$, and $\{b_m\}$ is a sequence positive numbers satisfying
\be\label{0.62}
\BBB_1=\sum_{m=1}^\iin \lm_m b_m^2<\iin,
\ee
where $\lm_m$ is the eigenvalue associated with $e_m$. It is well-known that under these assumptions, the Markov family corresponding to \ef{0.8} has a unique stationary measure $\mu^\es$ attracting the law of any solution with exponential rate (see the book \cite{KS-book}). We are interested in the asymptotics of the family $(\mu^\es)$ as the amplitude $\es$ of the noise goes to zero. Let us denote by $S(t)$ the semigroup acting on $H$ associated to \ef{0.8} with $\es=0$, and let $\ooo$ be the set of $\om$-limit points of $S(t)$. Thus, $u\in \ooo$ iff there is $u_0\in H$ and $t_k\to \iin$ such that $S(t_k)u_0\to u$ as $k\to\iin$. Note that this set is precompact in $H^2(D)\cap H^1_0(D)$ since it is a subset of the global attractor $\aaa$ of $S(t)$ which itself is compact in $H^2(D)\cap H^1_0(D)$ (see \cite{BV1992}).

\begin{definition}
We shall say that a set $E\subset H$ is \mpp{stochastically attracting} if $\mu^\es$ vanishes exponentially outside any neighborhood of $E$, that is we have\,\footnote{Note that if $(\mu^\es)$ satisfies the LDP, then a set is stochastically attracting iff its closure contains the kernel of the rate function.}
$$
\limsup_{\es\to 0}\es\ln\mu^{\es}(E^c_\eta)<0 \q\q \text{ for any }\eta>0,
$$
where $E_\eta$ stands for the open $\eta$-neighborhood of $E$ in $H$.
\end{definition}

What follows is the main result of this paper. 
\begin{theorem}\label{th1}
Under the above hypotheses, the family $(\mu^\es)$ is exponentially tight in $H^{1-\gamma}(D)$ for any $\gamma>0$. Moreover, the set $\ooo$ of $\om$-limit points of $S(t)$ is stochastically attracting. In particular, any weak limit of $(\mu^\es)$ is concentrated on the closure $\bar\ooo$. In the case of trivial limiting dynamics (i.e., when $\ooo$ is a singleton), the family $(\mu^\es)$ obeys the large deviations principle in $H$ (and by exponential tightness, also in $H^{1-\gamma}(D)\cap H$ for any $\gamma>0$).
\end{theorem}
Before outlining the main ideas behind the proof of this theorem, let us discuss the existing results in the literature on this subject. In the PDE setting, the study of this problem was initiated by Sowers in \cite{sowers-1992b}, where the author proves the LDP for stationary measures of the reaction-diffusion equation with a non-Gaussian perturbation. This was later extended to the case of multiplicative noise by Cerrai and R\"ockner in \cite{CeRo2005}. Recently, Brzezniak and Cerrai \cite{2015arXiv150900077B} established the LDP for the Navier-Stokes equation with additive noise on a 2D torus.  All these results cover only the case of trivial limiting dynamics, namely, the origin is the unique equilibrium of the deterministic equation (in which case the global attractor of a limiting equation is a singleton). The main ingredient here is the technique developed in \cite{sowers-1992b}. In the case when the limiting equation has arbitrary finite number of equilibria, we established the LDP for the nonlinear wave equation with smooth white noise in \cite{DM2015}. The proof relies on the development of Freidlin-Wentzell \cite{FW2012} and Khasminskii \cite{Khas2011} theories to the infinite-dimensional setting. Even though that equation is technically much more involved than the 2D Navier-Stokes equation, a crucial point in \cite{DM2015} is that due to the Lyapunov structure (i.e., existence of a function that decays on trajectories), we have an explicit description of the global attractor of the limiting equation. Namely, it consists of equilibrium points and joining them heteroclinic orbits.

\medskip
In the present work, for the first time, the aymptotics of stationary measures is studied, when the structure of the global attractor of the deterministic equation is not known. By adapting Sowers' argument, it is possible to show that stationary measures $(\mu^\es)$ vanish exponentially outside any neighborhood of the global attractor $\aaa$ of $S(t)$. Theorem \ref{th1} says that this decay holds true outside of a much smaller set, the set $\ooo$ of $\om$-limit points of $S(t)$. Note that even the simple consequence that any weak limit of $(\mu^\es)$ is concentrated on $\bar\ooo$ is new in the literature. 

\medskip
Let us describe in few words the main ideas of the proof. First, by developing Sowers' approach, we show that $(\mu^\es)$ satisfies the large deviations upper bound with a rate function $\vvv_\aaa$ (see \ef{0.30}) that has bounded level sets in $H^1(D)$ closed in $H$, and vanishes only on the global attractor $\aaa$. This immediately implies the exponential tightness in $H^{1-\gamma}(D)$ for any $\gamma>0$. Moreover, this also implies the LDP in the trivial case, when the set $\aaa$ is a singleton. Indeed, a simple argument shows that in this case $\vvv_\aaa$ provides also the large deviations lower bound. The main ingredients here are the classical Foia\c{s}-Prodi and exponential moment estimates for \ef{0.61}. We then use the exponential tightness in $H^{1-\gamma}(D)$ to prove the most involved part of our result, that $\ooo$ is stochastically attracting. To this end, we apply the mixing properties of the Markov process associated with \ef{0.8} to show that its stationary measure $\mu=\mu^\es$ can be recovered using only some local knowledge about $\mu$. Finally, let us mention that this reconstruction formula is proved in an abstract setting and, we think, might be of independent interest, see Section \ref{0.73}. 

\bigskip
{\bf Acknowledgments}. I thank Armen Shirikyan for very useful remarks.

\subsection*{Notation}
Given a Banach space $X$ and a positive constant $R$, we shall denote by $B_R(X)$ the closed ball of radius $R$ in $X$ centered at the origin. If $X=H$, we shall simply write $B_R$. For $u\in H$ and $\eta>0$, we write $B_\eta(u)$ for the open ball in $H$ of radius $\eta$ and centered at $u$. Similarly, for $A\subset H$, we shall denote by $A_\eta$ the open $\eta$-neighborhood of $A$ in $H$. We denote by $d(\cdot, \cdot)$ the distance in $L^2$, by $(\cdot, \cdot)$ its inner product, and by $\|\cdot\|$ the corresponding norm. Finally, we introduce
$$
H_{\vartheta}=\left\{v\in H: |v|_{H_\vartheta}^2=\sum_{j=1}^\iin b_j^{-2}\,(v,e_j)^2<\iin\right\}
$$
and $V=H^1_0(D)\cap H$. Note that in view of \ef{0.62}, we have $H_\vartheta\hookrightarrow V$ and the embedding is Hilbert-Schmidt (thus compact).

\section{Notions of exponential tightness and large deviations}
A family $(\mathfrak{m}^{\es})_{\es>0}$ of probability measures defined on a Polish space $X$ is said to be \mpp{exponentially tight} in $X$ if we have
\be\label{0.63}
\inf_{K}\limsup_{\es\to 0}\es \ln\mathfrak{m}^{\es}(X\backslash K)=-\iin,
\ee
where the infimum is taken over all compact sets $K$ in $X$. 

\medskip
A functional $\mathfrak{I}:X\to [0, \iin]$ is called \mpp{a (good) rate function} on $X$ if it has compact level sets, i.e., the set $\{\mathfrak{I}\le M\}$ is compact in $X$ for any $M\ge 0$. The family $(\mathfrak{m}^{\es})_{\es>0}$ obeys \mpp{large deviations principle} in $X$ with rate function $\mathfrak{I}$ if the following two properties hold:

\bi
\item{\it Lower bound}
\ei
$$
\liminf_{\es\to 0}\es\ln\mathfrak{m}^{\es}(G)\geq-\inf_{z\in G} \mathfrak{I}(z)\q \text{ for any } G\subset X \text{ open}.
$$

\bi
\item{\it Upper bound}
\ei
$$
\limsup_{\es\to 0}\es\ln \mathfrak{m}^{\es}(F)\leq-\inf_{z\in F} \mathfrak{I}(z)\q \text{ for any } F\subset X \text{ closed}.
$$
Note that a family satisfying large deviations upper bound is exponentially tight. An important property of exponentially tight family is that from its any sequence we can extract a subsequence that obeys the large deviations principle. As a corollary, if $(\mathfrak{m}^{\es})$ is an exponentially tight family of probability measures on $Y$ obeying large deviations in $X\hookleftarrow Y$, then it obeys large deviations in $Y$.

\medskip
In what follows, we shall say that $(\mathfrak{m}^{\es})$ is \mpp{weakly exponentially tight} in X if \ef{0.63} holds with the infimum taken over all bounded sets $K$ in $X$. If $(\mathfrak{m}^{\es})$ is weakly exponentially tight in a space $Y$ that is compactly embedded in $X$, then it is clear that $(\mathfrak{m}^{\es})$ is exponentially tight in $X$. 

\section{Proof of exponential tightness and LDP}
In this section, taking for granted some technical results established in the appendix, we shall prove that the family $(\mu^\es)$ of stationary measures is exponentially tight and in the case of trivial limiting dynamics obeys large deviations principle. Let us first introduce some notation. For $t\ge 0$, $v\in H$ and $\es>0$, we shall write $S^\es(t)v$ for the solution at time $t$ of equation \ef{0.8} issued from $v$. For $\ph\in L^2_{loc}(\rr_+;H)$, we shall denote by $S^\ph(t)v$ the solution at time $t$ of controlled equation
\be\label{0.26}
\dt u+Lu+B(u,u)=h+\ph
\ee
issued from $v$. For a trajectory $u_\cdot$ in $C(0, T;H)$, we introduce the energy
$$
I_T(u_\cdot)=J_T(\ph)=\f{1}{2}\int_0^T |\ph(s)|^2_{H_\vartheta}\dd s
$$
if there is $\ph\in L^2(0, T;H_\vartheta)$ such that $u_\cdot=S^\ph(\cdot)u_0$ and $I_T(u_\cdot)=\iin$ otherwise.

\subsection{Exponential tightness} To prove exponential tightness, we shall construct a function $\vvv_\aaa:H\to [0, \iin]$ with bounded level sets in $V$ and closed in $H$ that provides the large deviations upper bound for the family $(\mu^\es)$. Moreover, in the case of trivial limiting dynamics (and only then), this bound will imply the lower bound with the same function, and we get the LDP governed by (good) rate function $\vvv_\aaa$. Consider the semigroup $S(t):H\to H$ corresponding to 
\be\label{0.7}
\dt u+Lu+B(u,u)=h
\ee
and denote by $\aaa$ its global attractor. For $u_*\in H$, let $\vvv_\aaa(u_*)$ be the minimal energy needed to reach any neighborhood of $u_*$ from the set $\aaa$ in a finite time:
\be\label{0.30}
\vvv_\aaa(u_*)=\lim_{\eta\to 0}\inf\left\{I_s(u_\cdot), s>0, u_\cdot\in C(0, s;H): u_0\in\aaa, u_s\in B_\eta(u_*)\right\}.
\ee
Notice that the above limit (finite or infinite) exists, since the infimum written after the limit sign is monotone in $\eta>0$. We shall see below that this definition readily implies the closedness of level sets of $\vvv_\aaa$ in $H$.

\begin{proposition}\label{9.18}
Under the hypotheses of Theorem \ref{th1}, the function $\vvv_\aaa$ has bounded level sets in $V$ which are closed in $H$, and provides the large deviations upper bound for the family $(\mu^\es)$ in $H$, that is we have
\be\label{9.19}
\limsup_{\es\to 0}\es\ln \mu^{\es}(F)\leq-\inf_{u\in F} \vvv_\aaa(u)\q \text{ for any } F\subset H \text{ closed}.
\ee
In particular, the family $(\mu^\es)$ is weakly exponentially tight in $V$.
\end{proposition} 
By adapting Sowers' approach, it is easy to show (see Section \ref{0.20}) that bound \ef{9.19} will be established if we prove the following three properties:

\bi
\item \mpp{Trajectory inclusion}: for any positive constants $\De,\De'$ and $M$, there is $\eta>0$ such that
\be\label{0.12}
\{u(t): u(0)\in \aaa_\eta, I_t(u_\cdot)\leq M-\De'\}\subset K_{\De}(M), \q t> 0,
\ee
where $K_\De(M)$ is  the open $\De$-neighborhood of the level set $\{\vvv_\aaa\leq M\}$.

\item \mpp{Energy inequality}: for any positive constants $R$ and $\eta$, there is $T>0$ such that we have
\be\label{0.13}
a=\inf\{I_{T}(u_\cdot); \, u_\cdot\in C(0,T;H), \, u(0)\in B_R, \,u(T)\notin \aaa_\eta\}>0.
\ee

\item \mpp{Weak exponential tightness in $H$}: we have
\be\label{0.15}
\lim_{R\to\iin}\limsup_{\es\to 0}\es\ln\mu^\es(B_R^c)=-\iin.
\ee
\ei
Note that \ef{0.15} follows directly from Theorem 2.5.5 in \cite{KS-book} and Chebychev inequality.

\bigskip
{\it Compactness of level sets of $\vvv_\aaa$.} 

\medskip
{\it Step~1:} Let us first prove that the level sets $\{\vvv_\aaa\leq M\}$ are closed in $H$. 
To this end, let $u_*^j\in\{\vvv_\aaa\leq M\}$ be a sequence converging to $u_*$ in $H$ and let us show that $\vvv_\aaa(u_*)\leq M$. By definition of $\vvv_\aaa$, we need to prove that for any positive constants $\eta$ and $\eta'$ there is an initial point $u_0\in\aaa$, a finite time $T>0$, and an action function $\ph$ such that
\be\label{0.78}
J_{T}(\ph)\leq M+\eta'\q \text{ and }\q \|S^{\ph}(T)u_{0}-u_*\|\leq \eta.
\ee
We fix $j$ so large that
\be\label{0.79}
\|u_{*}^j-u_*\|\leq \eta/2.
\ee

Since $\vvv_\aaa(\uu_*^j)\leq M$, there exist a point $u_0\in\aaa$, a time $T>0$ and an action $\ph$ such that
$$
J_{T}(\ph)\leq M+\eta'\q \text{ and }\q \|S^{\ph}(T)u_{0}-u_*^j\|\leq \eta/2.
$$
Combining this with inequality \eqref{0.79}, we infer \eqref{0.78}.

\medskip
{\it Step~2:} We now show that the levels sets are bounded in $V$, that is, for any $M\ge 0$ there is $R(M)>0$ such that

\be\label{0.52}
\{\vvv_\aaa\le M\}\subset B_{R(M)}(V).
\ee
To this end, let us fix $u_*\in \{\vvv_\aaa\le M\}$. By definition of $\vvv_\aaa$, there is a sequence of initial points $u_0^j\in\aaa$, of positive times $s_j$, and functions $\ph^j\in L^2(0, s_j; H_{\vartheta})$ with energy
$$
\int_0^{s_j}|\ph^j(\tau)|_{H_\vartheta}^2\dd \tau\le M+1
$$ 
such that 
$$
d(S^{\ph^j}(s_j)u_0^j, u_*)\to 0\q\text{ as } j\to\iin.
$$
In view of Lemma \ref{0.51}, there is a positive constant $R(M)$ such that
$$
\sup_{t\in [0, s_j]}|S^{\ph^j}(t)u_0^j|_{V}\le R(M)\q\text{ for any } j\ge 1.
$$
Taking $t=s_j$ in this inequality and using the above convergence, we see that the $V$-norm of $u_*$ is bounded by $R(M)$ and thus infer \ef{0.52}.

\bigskip
{\it Proof of the trajectory inclusion}: {\it Step~1:} We note that this inclusion is trivial for $\eta=0$. To prove that it holds also for $\eta>0$ sufficiently small, we shall apply the Foia\c{s}-Prodi estimate for the Navier-Stokes equation. Assume that \ef{0.12} is not true, and let us find sequences of positive numbers $T_m$ and $\eta_m\to 0$, of initial points $u^m_0\in \aaa_{\eta_m}$ and of action functions $\ph^m$ with $J_{T_m}(\ph^m)\le M-\De'/2$ such that the flow $u^m(t)=S^{\ph^m}(t)u_0^m$ satisfies 
\be\label{0.25}
u^m(T_m)\notin K_\De(M).
\ee 
It is easy to see that
\be\label{0.23}
\int_0^{t}\|\g u^m(s)\|^2\dd s\le \mmm(1+t) \q\text{ for } t\in [0, T_m],
\ee
where the constant $\mmm>0$ depends only on $\|h\|$ and $M$.
For any $m\ge 1$, let us fix $w_0^m\in\aaa\cap B_{\eta_m}(u_0^m)$ and introduce an intermediate flow $w^m(t)$ defined on the time interval $[0, T_m]$ that solves 
\be\label{0.21}
\dt w^m+Lw^m+B(w^m,w^m)=h+\ph^m+\lm P_N(u^m-w^m), \q w^m(0)=w^m_0,
\ee
where $\lm>0$ and $N\in\nn$ are some constants, and $P_N$ stands for the orthogonal projection from $H$ to its subspace spanned by the first $N$ eigenfunctions of the Stokes operator. Let us use inequality \ef{0.23} together with Theorem 2.1.28 of \cite{KS-book}, to choose $\lm=\lm(\mmm)$ and $N=N(\mmm)$ such that on the interval $[0, T_m]$, we have
\be\label{0.22}
\|u^m(t)-w^m(t)\|^2\le e^{-t+c\mmm}\|u^m_0-w^m_0\|^2\le e^{-t+c\mmm}\eta_m^2,
\ee
with an absolute constant $c>0$.

\medskip
{\it Step~2:} Let us estimate $I_{T_m}(w_\cdot^m)$. By the very definition of $I_T$, we have
\begin{align}\label{0.24}
I_{T_m}(w^m_\cdot)&=\f{1}{2}\int_0^{T_m}|\ph^m(s)+\lm P_N\left(u^m(s)-w^m(s)\right)|_{H_\vartheta}^2\dd s\notag\\
&\le\f{a}{2}\int_0^{T_m}|\ph^m(s)|_{H_\vartheta}^2\dd s+\f{a}{2(a-1)}\int_0^{T_m}|\lm P_N\left(u^m(s)-w^m(s)\right)|_{H_\vartheta}^2\dd s\notag\\
&\le a J_{T_m}(\ph^m)+\f{a}{2(a-1)}\lm^2 C(N)\int_0^{T_m}\|u^m(s)-w^m(s)\|^2\dd s\notag,
\end{align}
where $a>1$ is any constant. Using this with inequality \ef{0.22} together with the fact that the constants $\lm$ and $N$ depend only on $\mmm$, we get
$$
I_{T_m}(w^m_\cdot)\le aJ_{T_m}(\ph^m)+\f{a}{a-1}C(\mmm)\eta_m^2.
$$

\medskip
{\it Step~3:}
Choosing $a=(M-\De'/4)/(M-\De'/2)$, we derive $I_{T_m}(w^m_\cdot)\le M$ for all $m$ sufficiently large. Since $w^m_0\in\aaa$, we obtain that the point $w^m(T_m)$ is reached from the global attractor $\aaa$ at finite time with energy not bigger than $M$. By definition of $\vvv_\aaa$, this implies $w^m(T_m)\in \{\vvv_\aaa\le M\}$. Combining this with inequalities \ef{0.25} and \ef{0.22}, we arrive at a contradiction. Inclusion \ef{0.12} is thus established.

\bigskip
{\it Derivation of energy inequality}: As above, we proceed by contradiction. If inequality \ef{0.13} is not true, then there are positive constants $R$ and $\eta$ such that 
$$
\inf\{I_{m}(u_\cdot); \, u_\cdot\in C(0,m;H), \, u(0)\in B_R, \,u(m)\notin \aaa_\eta\}=0, \q m\in\nn.
$$
For each $m\geq 1$, let us find $u^m_0\in B_R$ and action $\ph^m$ defined on the time interval $[0,m]$ with energy $J_m(\ph^m)$ smaller than $e^{-m^2}$ such that the flow $u^m(t)=S^{\ph^m}(t)u^m_0$ satisfies 
\be\label{2.4}
u^m(m)\notin \aaa_\eta.
\ee
Let us set $v^m(t)=S(t)u^m_0$. In view of Lemma \ref{0.16}, we have
\be\label{0.18}
\|u^m(t)-v^m(t)\|^2\le C_Re^{ct\|h\|^2}\int_0^t \|\ph(s)\|^2\dd s\le C_R' e^{ct\|h\|^2}J_t(\ph^m).
\ee
Taking $t=m$ in this inequality and using $J_m(\ph^m)\le e^{-m^2}$, we see that the distance $\|u^m(m)-v^m(m)\|$ converges to zero as $m$ goes to infinity. Combining this with \ef{2.4}, we get
$$
v^m(m)\notin \aaa_{\eta/2}
$$
for all $m$ sufficiently large. However, since $\aaa$ is the global attractor of the semigroup $S(t)$, we have
$$
d(v^m(m), \aaa)\le \sup_{u_0\in B_R}d(S(m)u_0, \aaa)\to 0\q\text{ as }m\to\iin.
$$ 
Inequality \ef{0.13} is proved.

\subsection{The LDP in the case of trivial limiting dynamics}\label{0.59}
Here we show that in the case when the global attractor is a singleton, the function $\vvv_\aaa$ given by \ef{0.30} provides also a lower bound for $(\mu^\es)$ and thus governs the LDP of that family.

\bigskip
Let $\aaa=\{\hat u\}$. In view of Proposition \ref{9.18}, the family $(\mu^\es)$ is tight. Moreover, since function $\vvv_\aaa$ vanishes only on $\aaa$, any weak limit of this family is concentrated on $\aaa=\{\hat u\}$. It follows that
\be\label{0.57}
\mu^\es\rightharpoonup\De_{\hat u}.
\ee 
In order to get the lower bound, it is sufficient to prove that for any $u\in H$ and any positive constants $\De$ and $\De'$, there is $\es_*>0$ such that we have 
\be\label{0.58}
\mu^\es(B_\De(u))\geq \exp(-(\vvv_\aaa(u)+\De')/\es) \q\text {for }\es\leq \es_*.
\ee
We may assume that $\vvv_\aaa(u)<\iin$. By definition of functional $\vvv_\aaa$, there is a time $s>0$ and function $\ph\in L^2(0,s;H_\vartheta)$ such that
$$
J_s(\ph)\leq \vvv_{\aaa}(u)+\De'\q\text{ and }\q S^\ph(s)\hat u\in B_{\De/4}(u).
$$
Due to continuity of $S^\ph$ with respect to the initial point, there is $\eta>0$ such that for $v\in B_{\eta}(\hat u)$, we have $S^\ph(s)v\in B_{\De/2}(u)$. 
Now using the stationarity of $\mu^\es$ and Theorem \ref{0.54}, we derive
\begin{align*}
\mu^\es(B_\De(u))&=\int_H \pp\left\{S^\es(s)v\in B_\De(u)\right\}\mu^\es(\Dd v)\ge \int_{B_\eta(\hat u)} \pp\left\{S^\es(s)v\in B_\De(u)\right\}\mu^\es(\Dd v)\\
&\ge\mu^\es(B_\eta(\hat u)) \exp(-(\vvv_{\aaa }(u)+2\De')/\es).
\end{align*}
Combining this inequality with convergence \ef{0.57} and using the portmanteau theorem, we get
$$
\mu^\es(B_\De(u))\geq C(\eta) \exp(-(\vvv_{\aaa}(u)+2\De')/\es)\geq  \exp(-(\vvv_{\aaa}(u)+3\De')/\es)
$$
for $\es$ sufficiently small.
Since $\De'$ was arbitrary, this is equivalent to \ef{0.58}.

\section{Abstract result}\label{0.73}
Here we prove a formula for recovering a stationary measure of a Markov process with some good mixing properties. This will be used in the next section to establish stochastic attractiveness of $\ooo$.
We first introduce some notation and terminology. Given a metric space $X$, we shall denote by $\bB_0(X)$ the space of bounded measurable functions on $X$ endowed with the following convergence: we shall say that a sequence $\psi_n$ converges to $\psi$ in $\bB_0(X)$ (or that $\psi_n$ $\bB$-converges to $\psi$) if 
$$
\sup_n\sup_{v\in X}|\psi_n(v)|<\iin
$$
and 
$$
\sup_{v\in B} |\psi_n(v)-\psi(v)|\to 0\q\text{ as }n\to\iin
$$
for any bounded set $B\subset X$.

\medskip
Let $(u_t, \pp_v)_{t\in \rr_+}$ be a Markov process in a metric space $X$ possessing an invariant measure  $\mu\in \ppp(X)$. 
We shall say that $\mu$ is \mpp{mixing in} $\bB_0(X)$ if for any bounded Lipschitz continuous function $\psi:X\to \rr$, there is $t_n\to\iin$ such that the sequence $(P_{t_n}\psi)$ converges to $(\psi, \mu)$ in $\bB_0(X)$, where $P_t$ stands for the corresponding Markov operator. Note that if $\mu$ is mixing in $\bB_0(X)$ for $(u_t, \pp_v)$, then it is the unique invariant measure of that process.

\medskip
Let us be given a Markov process $(u_t,\pp_v)_{t\ge 0}$ and a family of $\fff_t$-stopping times $\{\tau(v)\}_{v\in X}$, where $\fff_t$ is the filtration generated by $u_t$. We shall say that \mpp{condition} ({\bf A}) is fulfilled for $(u_t,\pp_v)_{t\ge 0}$ and $\{\tau(v)\}_{v\in X}$ if we have the following.

\bi

\item The process $(u_t, \pp_v)_{t\in \rr_+}$ is defined on a Polish space $X$ and possesses an invariant measure $\mu\in \ppp(X)$ that is mixing in $\bB_0(X)$. Moreover, we assume that for any compact set $K\subset X$, any $T>0$ and $\eta>0$ there is a bounded set $B\subset X$ such that
\be\label{0.74}
\pp_v\left\{(u_t)_{t\in [0, T]}\subset B \right\}\ge 1-\eta \q\text{ for any } v\in K.
\ee

\item 
The family $\{\tau(v)\}_{v\in X}$ satisfies
\be\label{0.75}
\lim_{s\to\iin}\sup_{v\in K}\pp_v\left\{\tau(v)\ge s\right\}=0
\ee
for any compact set $K\subset X$.
Moreover, the map $\tau:(\om, v)\to \tau^\om(v)$ is measurable from the product space $\Omega\times X$ to $[0, \iin]$.
 \ei

\bigskip
\nt
For any $\De>0$, $\psi\in \bB_0(X)$ and $v\in X$, introduce the operator
\be\label{0.65}
\R_\De\psi(v)=\f{1}{\De}\e_v\int_{\tau}^{\tau+\De} \psi(u_t)\dd t\equiv\f{1}{\De}\e_v\int_{\tau(v)}^{\tau(v)+\De} \psi(u_t)\dd t.
\ee
Note that thanks to the measurability of $\tau$, $\R_\De\psi$ is measurable from $X$ to $\rr$. Indeed, we have
$$
\R_\De\psi(v)=\f{1}{\De}\e\int_{\rr_+}\ch_{[\tau(v), \tau(v)+\De]} \psi(u_t(v))\dd t,
$$
where $u_t(v)$ stands for the trajectory issued from $v$. The Fubini-Tonelli theorem allows to conclude.

\medskip
The following proposition is the main result of this section.

\begin{proposition}\label{0.66}
Let condition \text{\textnormal{({\bf A})}} be fulfilled for a Markov process $(u_t, \pp_v)_{t\in \rr_+}$ and family $\{\tau(v)\}_{v\in X}$. Then, for any positive constant $\De$, the map $\lm$ given by relation
\be\label{0.64}
\lm(\psi)=(\R_\De\psi, \mu)\equiv \int_X \R_\De\psi(v)\mu(\Dd v),
\ee
is continuous from $\bB_0(X)$ to $\rr$ and satisfies
\be\label{0.69}
\mu(\dt \Gamma)\le \lm(\dt\Gamma)\le \lm(\bar\Gamma)\le\mu(\bar\Gamma) \q\text{ for any }\Gamma\subset X,
\ee
where $\dt\Gamma$ and $\bar\Gamma$ stand for its interior and closure, respectively, and we write $\lm(\Gamma)$ for the value $\lm(\ch_\Gamma)$.

\end{proposition}
\bp
\nt

\medskip
{\it Step~1:} To simplify the notation, we shall assume that $\De=1$. We first use the Khasminskii type argument (see Chapter 4 in \cite{Khas2011}). For any $\psi\in\bB_0(X)$, $v\in X$ and $s>0$, we have
\begin{align*}
\e_v\int_\tau^{\tau+1}\psi(u_{t+s})\dd t&=\e_v\int_0^\iin \ch_{[\tau, \tau+1]}(t)\psi(u_{t+s})\dd t\\
&=\int_0^\iin\e_v[\ch_{[\tau, \tau+1]}(t)\psi(u_{t+s})]\dd t,
\end{align*}
where we write $\tau$ for $\tau(v)$.
Moreover, since $\tau$ is an $\fff_t$-stopping time, the characteristic $\ch_{[\tau, \tau+1]}(t)$ is $\fff_t$-measurable for any $t\ge 0$. It follows that
$$
\e_v[\ch_{[\tau, \tau+1]}(t)\psi(u_{t+s})]=\e_v\e_v[\ch_{[\tau, \tau+1]}(t)\psi(u_{t+s})|\fff_t]=\e_v[\ch_{[\tau, \tau+1]}(t)P_s\psi(u_t)],
$$
whence we infer
\be\label{0.67}
\e_v\int_\tau^{\tau+1}\psi(u_{t+s})\dd t=\int_0^\iin \e_v[\ch_{[\tau, \tau+1]}(t)P_s\psi(u_t)]\dd t=\e_v\int_\tau^{\tau+1}P_s\psi(u_t)\dd t.
\ee

\medskip
{\it Step~2:} Now let $\lm:\bB_0(X)\to\rr$ be given by
$$
\lm(\psi)=(\R_1\psi, \mu)=\left(\e_\cdot\int_\tau^{\tau+1}\psi(u_t)\dd t, \mu\right).
$$

\medskip
Thanks to \ef{0.67}, we have
\begin{align*}
\lm(P_s\psi)&=(\R_1 P_s\psi, \mu)=\left(\e_\cdot\int_\tau^{\tau+1}P_s\psi(u_t)\dd t, \mu\right)\\
&=\left(\e_\cdot\int_\tau^{\tau+1}\psi(u_{t+s})\dd t, \mu\right)=\left(\e_\cdot\int_{\tau+s}^{\tau+1+s}\psi(u_{t})\dd t, \mu\right)\\
&=\lm(\psi)+\left(\e_\cdot\int_{\tau+1}^{\tau+1+s}\psi(u_{t})\dd t, \mu\right)-\left(\e_\cdot\int_{\tau}^{\tau+s}\psi(u_{t})\dd t, \mu\right).
\end{align*}
Conditioning with respect to $\fff_1$ and using the stationarity of $\mu$, we see that the last two terms are equal, so that
\be\label{0.68}
\lm(P_s\psi)=\lm(\psi)
\ee
for any $s>0$ and $\psi\in \bB_0(X)$.

\medskip
{\it Step~3:} We now prove that $\lm$ is continuous from $\bB_0(X)$ to $\rr$. Let $(\psi_n)$ be a sequence $\bB$-converging to zero and let us show that $\lm(\psi_n)\to 0$. We may assume that $(\psi_n)$ is uniformly bounded by $1$. Let us fix any $\eta>0$. Since $X$ is Polish, we can use Ulam's theorem, to find $K\subset X$ compact such that $\mu(K^c)\le\eta$. It follows that
$$
|\lm(\psi_n)|\le\left(\ch_K(\cdot)\e_\cdot\int_\tau^{\tau+1}|\psi_n(u_t)|\dd t, \mu\right)+\eta.
$$
Now let us use \ef{0.75} to find $R>0$ such that
$$
\sup_{v\in K}\pp_v\left\{\tau\ge R\right\}\le \eta.
$$
Once $R$ is fixed, we use \ef{0.74} to find a bounded set $B\subset X$ such that
$$
\pp_v\left\{(u_t)_{t\in [0, R+1]}\subset B\right\}\ge 1-\f{\eta}{R+1} \q\text{ for }v\in K.
$$
It follows that for any $v\in K$, we have
\begin{align*}
\e_v\int_\tau^{\tau+1}|\psi_n(u_t)|\dd t&\le \pp_v\left\{\tau\ge R\right\}+\e_v\left(\ch_{\tau< R}\int_\tau^{\tau+1}|\psi_n(u_t)|\dd t\right)\\
&\le \eta+\e_v\left(\int_0^{R+1}|\psi_n(u_t)|\dd t\right)\le 2\eta+ (R+1)\sup_{B}|\psi_n|.
\end{align*}
We thus derive
$$
|\lm(\psi_n)|\le 3\eta+(R+1)\sup_{B}|\psi_n|.
$$
Since $B$ is bounded and $\psi_n$ $\bB$-converges to zero, the second summand in this inequality is smaller than $\eta$ for $n$ sufficiently large. Now recalling that $\eta$ was arbitrary, we infer that $\lm(\psi_n)\to 0$.

\medskip
{\it Step~4:} Let us fix a bounded Lipschitz continuous function $\psi:X\to\rr$. Since $\mu$ is mixing in $\bB_0(X)$, the sequence $(P_{t_n}\psi)$ converges to $(\psi, \mu)$ in $\bB_0(X)$ for some $t_n\to\iin$. By continuity of $\lm$ from $\bB_0(X)$ to $\rr$ and \ef{0.68}, we get 
$$
\lm(\psi)=\lm(P_{t_n}\psi)\to \lm((\psi, \mu))=(\psi, \mu)\lm(1)=(\psi, \mu).
$$
Now fixing a closed subset $F\subset X$, approximating the characteristic function of $F$ by Lipschitz continuous functions $\ch_F\le \psi_n\le 1$ in the sense of pointwise convergence, and using the Lebesgue theorem on dominated convergence, we see that $\lm(F)\le \mu(F)$, which implies \ef{0.69}.
\ep

\section{The set of $\om$-limit points is stochastically attracting}
Here we prove that $(\mu^\es)$ decays exponentially outside any neighborhood of the set $\ooo$ of $\om$-limit points of $S(t)$, that is
\be\label{0.76}
\limsup_{\es\to 0}\es\ln\mu^{\es}(\ooo^c_\eta)<0 \q\q \text{ for any }\eta>0.
\ee
{\it Proof of \ef{0.76}.} We shall derive this result from Proposition \ref{0.66}. 

\medskip
{\it Step~1:} Let us fix  any $\eta>0$ and for $v\in H$, denote by $\tau(v)$ the first instant when the deterministic flow $S(t)v$ hits the set $\ooo_{\eta/4}$. We claim that the process $u_t=S^\es(t)$ and  family $\{\tau(v)\}_{v\in H}$ satisfy \mpp{condition} ({\bf A}) of previous section. Indeed, since $\tau(v)$ is constant for any $v\in H$, the set $\{\tau(v)\le t\}$ is either empty or is the whole probability space $\Omega$, so $\tau(v)$ is adapted to any filtration. By the same reason, to show that $\tau$ is measurable on the product $\Omega\times H$, it is sufficient to prove that it is measurable from $H$ to $\rr$. The latter follows from the upper semi-continuity of $\tau$. Indeed, the set $\left\{v\in H: \tau(v)< a\right\}$ is open in $H$ for any $a>0$, since
\begin{align*}
\left\{v\in H: \tau(v)< a\right\}&=\left\{v\in H: \exists\, t< a, S(t)v\in \ooo_{\eta/4}\right\}\\
&=\bigcup_{t<a } \left\{v\in H: S(t)v\in \ooo_{\eta/4}\right\}\\
&=\bigcup_{t<a, t\in \qq} \left\{v\in H: S(t)v\in\ooo_{\eta/4}\right\},
\end{align*}
where we used the continuity of $S(t)$.
Thanks to \ef{0.9}, relation \ef{0.75} is also fulfilled.

\medskip
Further, the process $u_t=S^\es(t)$ satisfies \ef{0.74} in view of the supermartingale type inequality (see Proposition 2.4.10 in \cite{KS-book}). Moreover, the corresponding invariant measure $\mu=\mu^\es$ is mixing in $\bB_0(H)$. Indeed, by Theorem 3.5.2 of \cite{KS-book}, there are positive constants $C$ and $\al$ such that for any 1-Lipschitz continuous function $\psi:H\to \rr$, we have
$$
|P_t\psi(v)-(\psi,\mu)|\le Ce^{-\al t}\left(1+\|v\|^2\right) \q \text{ for any } t\ge 0, v\in H.
$$
In particular, for any $t_n\to\iin$, the sequence $(P_{t_n}\psi)$ converges to $(\psi, \mu)$ uniformly on any bounded subset of $H$. Combining this with the fact that $(P_{t_n}\psi)$ is uniformly bounded in $H$ (by $|\psi|_\iin$), we infer that $(P_{t_n}\psi)\to (\psi,\mu)$ in $\bB_0(H)$. Thanks to Proposition \ref{0.66}, for any $\De>0$, we have

$$
\mu^\es(\bar\ooo_\eta^c)\le\lm^\es(\bar\ooo_\eta^c),
$$
where 
$$
\lm^\es(\Gamma)=\f{1}{\De}\left(\e_\cdot\int_{\tau}^{\tau+\De}\ch_{\Gamma}(S^\es(t))\dd t, \mu^\es\right).
$$

\medskip
{\it Step~2:}
Let us use weak exponential tightness of $(\mu^\es)$ in $V$, to find $R>0$ such that for $K=B_R(V)$, we have
\be\label{0.80}
\limsup_{\es\to 0}\es\ln\mu^{\es}(K^c)<0.
\ee
Clearly, $K$ is compact in $H$ and, by the previous step, we have
\be\label{0.81}
\mu^\es(\bar\ooo_\eta^c)\le \f{1}{\De}\left(\ch_{K}(\cdot)\e_\cdot\int_{\tau}^{\tau+\De}\ch_{\bar\ooo_\eta^c}(S^\es(t))\dd t, \mu^\es\right)+\mu^\es(K^c).
\ee
 Introduce the time
$$
T=\sup_{v\in K}\tau(v)+1
$$
and the event
$$
A_v=\left\{\om\in \Omega: d_{C(0, T;H)}(S^\es(t)v, \{I_T\le a\})\le \eta/4\right\},
$$
where $a>0$.
Thanks to Theorem \ref{0.54}, for $\es>0$ sufficiently small, we have
\be\label{0.82}
\sup_{v\in K}\pp(A_v^c)\le \exp(-a/2\es).
\ee
Let us find $a=a(\eta, T,K)>0$ so small that for any curve $u_\cdot\in \{I_T\le a\}$ in $C(0, T;H)$ issued from a point $v\in K$, we have $d_{C(0, T; H)}(u_\cdot, S(t)v)\le \eta/4$. Clearly, such choice is possible. It follows that for any $v\in K$, on the event $A_v$, we have
$$
d_{C(0, T; H)}(S^\es(t)v, S(t)v)\le \eta/2.
$$

\medskip
{\it Step~3:}
Here we show that there is $\De>0$ so small, that for any curve $u_\cdot$ issued from a point $v\in K$ and lying in the $\eta/2$-neighborhood of $S(t)v$ in the space $C(0, T; H)$, we have $u_t\in\bar\ooo_{\eta}$ for $t\in [\tau(v), \tau(v)+\De]$. To this end, it is sufficient to prove that there is $\De>0$ such that $S(t)v\in \bar\ooo_{\eta/2}$ for any $v\in K$ and $t\in[\tau(v), \tau(v)+\De]$. Assume the opposite and find $\De_j\to 0$ and $v_j\in K$ such that 
\be\label{0.77}
S(\tau(v_j)+\De_j)v_j\notin\bar\ooo_{\eta/2}, \q j\ge 1.
\ee
Note that $S(\tau(v_j)+\De_j)v_j=S(\De_j)w_j$, where $w_j=S(\tau(v_j))v_j\in \bar\ooo_{\eta/4}$. Moreover, since $v_j\in K\equiv B_R(V)$, there is $C_R>0$ such that
$$
|w_j|_V\le\sup_{t\in [0, T]}|S(t)v_j|_V\le C_R \q \text{ for }j\ge 1.
$$
Combining this with the compactness of the embedding $V\hookrightarrow H$, we may assume that $w_j$ converges to some $w_*$ in $H$. In particular, we have $w_*\in \bar\ooo_{\eta/4}$ and 
$$
d_{C(0, 1;H)}(S(t)w_j, S(t)w_*)\to 0.
$$
By the triangle inequality, we get
\begin{align*}
d(S(\tau(v_j)+\De_j)v_j, \bar\ooo_{\eta/3})&\le d(S(\De_j)w_j, S(\De_j)w_*)+d(S(\De_j)w_*, \bar\ooo_{\eta/3})\\
&\le d_{C(0, 1;H)}(S(t)w_j, S(t)w_*)+d(S(\De_j)w_*, \bar\ooo_{\eta/3})\to 0.
\end{align*}
This clearly contradicts \ef{0.77} and proves our assertion concerning the existence of such $\De>0$.

\medskip
{\it Step~4:}
 It follows from the previous two steps that for $v\in K$, on the event $A_v$, we have 
 $$
 S^\es(t)v\in \bar\ooo_{\eta} \q\text{ for } t\in [\tau(v), \tau(v)+\De].
 $$
 Therefore, for such choice of $\De$, the quantity 
 $$
\ch_{K}(v)\e_v\int_{\tau}^{\tau+\De}\ch_{\bar\ooo_\eta^c}(S^\es(t))\dd t
 $$
 vanishes on $A_v$. Using this together with \ef{0.81}-\ef{0.82}, we infer
 $$
 \mu^\es(\bar\ooo_\eta^c)\le \exp(-a/2\es)+\mu^\es(K^c).
 $$
Finally, combining last inequality with \ef{0.80}, we derive
 $$
\limsup_{\es\to 0}\es\ln\mu^\es(\bar\ooo_\eta^c)<0,
$$
and since $\eta>0$ was arbitrary, this is equivalent to \ef{0.76}.

\subsection{Additional bound for hitting times}
The goal of this section is to establish another estimate for hitting times. We denote by $\tau^\es_\eta(v)$ the first instant when the trajectory $S^\es(t)v$ hits the set $\bar\ooo_\eta$.
\begin{lemma}
For any $\eta>0$ and $R>0$, we have 
\be\label{0.72}
\lim_{s\to\iin}\limsup_{\es\to 0}\sup_{v\in B_R}\es\ln\pp\left(\tau^\es_\eta(v)\ge s\right)=-\iin.
\ee
\end{lemma}
\bp
\nt

\medskip
{\it Step~1: Reduction.}
As is shown in the derivation of Lemma 2.3 in \cite{DM2015}, using the Markov property and supermartingale inequality, the proof of \ef{0.72} can be reduced to 
$$
\lim_{s\to\iin}\limsup_{\es\to 0}\sup_{v\in B_R}\es\ln\pp\left(\tau^\es_\eta(v)\ge s\right)<0.
$$
On the other hand, thanks to large deviations for trajectories, it is sufficient to prove that
\be\label{0.9}
\sup_{v\in B_R}\l_\eta(v)<\iin,
\ee
where $l_v(\eta)$ is the first instant when $S(t)v$ hits the set $\bar\ooo_\eta$.

\medskip
{\it Step~2: Derivation of inequality \ef{0.9}.}
Assume the opposite and let us find $R>0$ and $\eta>0$ for which this inequality fails. 
Then, there exists a sequence $(v_m)\subset B_{R}$ such that 
\be\label{8.39}
l_\eta(v_m)\geq 2m.
\ee
Since $\aaa$ is absorbing for $S(t)$, we have
$$
d_m=d(S(m)v_m, \aaa)\to 0 \q\text{ as }m\to\iin. 
$$
Let us find $w_m\in \aaa$, such that
$$
d(S(m)v_m, w_m)=d_m.
$$
Further, since $\aaa$ is compact, there is $(m_k)\subset\nn$ and $w_*\in \aaa$, such that $d(w_{m_k}, w_*)\to 0$. Combining this with the triangle inequality, we get
$$
d(S(m_k)v_{m_k}, w_*)\le d_{m_k}+d(w_{m_k}, w_*)\to 0 \q\text{ as }k\to\iin.
$$
By the continuity of $S(t)$, for any $s>0$, we have
$$
S(s+m_k)v_{m_k}=S(s)S(m_k)v_{m_k}\to S(s) w_*.
$$
We fix $s>0$ so large that
$$
S(s)w_*\in \ooo_{\eta/2}.
$$
In view of \eqref{8.39}, $S(s+m_k)v_{m_k}\notin \bar\ooo_\eta$ for $k\geq 1$ large enough. This contradicts the above two relations and proves inequality \ef{0.9}.
\ep

\section{Appendix}
In this section, we collected some technical results used in the main text. 
\subsection{Large deviations for trajectories}
Let us fix a closed bounded set $B$ in $H$ and let $T>0$. Introduce the Banach space $\mathsf{\yyy}_{B,T}$ of continuous functions $y(\cdot,\cdot):B\times[0,T]\to H$ endowed with the topology of uniform convergence. The following result is classical, see for instance the book \cite{DZ1992} and paper \cite{CM-2010} (see also the remark on the uniformity after Theorem 6.2 in \cite{DM2015}).

\begin{theorem}\label{0.54}
Under the hypotheses of Theorem \ref{th1}, $(S^{\es}(t)v, t\in [0, T], v\in B)_{\es>0}$ regarded as a family of random variables in $\yyy_{B,T}$ satisfies the large deviations principle with rate function $I_T:\yyy_{B,T}\to [0,\iin]$ given by 
$$
I_T(y(\cdot,\cdot))=\f{1}{2}\int_0^T |\ph(s)|_{H_\vartheta}^2\dd s
$$
if there is $\ph\in L^2(0,T;H_\vartheta)$ such that $y(t, v)=S^{\ph}(t)v$, and equal to $\iin$ otherwise.
\end{theorem}
\subsection{Derivation of the upper bound \ef{9.19} using \ef{0.12}-\ef{0.15}}\label{0.20}
We use the argument developed in \cite{DM2015} that relies on some ideas introduced by Sowers in \cite{sowers-1992b}. It is well-known (e.g., see Chapter 12 of \cite{DZ1992}) that \ef{9.19} is equivalent to the following. For any positive constants $\De, \De'$ and $M$ there is $\es_*>0$ such that
\be\label{0.55}
\mu^\es(u\in H: d(u,\{\vvv_\aaa\leq M\})\geq\De)\leq\exp(-(M-\De')/\es)\q\text{ for } \es\leq\es_*.
\ee
The constants $\De, \De'$ and $M$ are assumed to be fixed, from now on. 

\bigskip
{\it Step~1:} Let us find $\eta>0$ such that we have inclusion \ef{0.12}, where $\De$ should be replaced by $\De/2$. Also, let us use \ef{0.15}, to find a positive constant $R$ such that 
\be\label{0.53}
\iI_1:=\mu^\es(B_R^c)\le \exp(-M/\es)
\ee
for $\es>0$ sufficiently small.
Once these constants are fixed, we find $T>0$ such that we have \ef{0.13}.
For any $m\geq 1$, we introduce the set 
$$
\ees_m=\{u_\cdot\in C(0, mT;H): u(0)\in B_R; \,\,u(jT)\in  \aaa_\eta^c\cap B_R, \, j=1,\ldots, m\}.
$$
The structure of this set and inequality \ef{0.13} imply that
$$
\inf\{I_{mT}(u_\cdot); u_\cdot\in \ees_m\}>M
$$
for $m>(M+1)/a$. It follows from Theorem \ef{0.54} that 
\be
\iI_2:=\sup_{v\in B_R}\pp\left\{S^\es(\cdot)v\in \ees_m\right\}\le \exp(-M/\es).
\ee

\bigskip
{\it Step~2:}
Let us show that for $t_1=(m+1)T$, we have
\be
\iI_3:=\int_{B_R}\pp\left\{S^{\es}(t_1)v\notin K_{\De}(M),\,S^{\es}(\cdot)v\notin \ees_m\right\}\mu^\es(\Dd v)\leq \exp(-(M-2\De')/\es).\label{0.56}
\ee
Indeed, we have
\begin{align*}
\iI_3
&\leq \sum_{k=1}^m \int_{B_R}\pp\left\{S^{\es}(t_1)v\notin K_{\De}(M),\,S^{\es}(kT)v\in\aaa_\eta\right\}\mu^\es(\Dd v)\notag\\
&\q+ \sum_{k=1}^m \int_{B_R}\pp\left\{S^{\es}(t_1)v\notin K_{\De}(M),\,S^{\es}(kT)v\in B_R^c \right\}\mu^\es(\Dd v)=\iI_3'+\iI_3''.
\end{align*}
It follows from inclusion \ef{0.12} that
$\iI_3'\leq\exp(-(M-\De')/\es)$.
On the other hand, thanks to stationarity of $\mu^\es$, we have
$$
\iI_3''\leq\sum_{k=1}^m\int_{H}\pp\left\{S^{\es}(kT)v\in B_R^c  \right\}\mu^\es(\Dd v)=m\,\iI_1\le m\exp(-M/\es),
$$
which leads to \ef{0.56}. 

\bigskip
{\it Step~3:} We are ready to derive \ef{0.55}. Indeed, it follows from definition of $K_\De(M)$ and stationarity of $\mu^\es$ that
\begin{align*}
\mu^\es(u\in H: d(u,\{\vvv_\aaa\leq M\})\ge\De)&=\mu^\es(u\in H: u\notin K_\De(M))\\
&=\int_{H}\pp\left\{S^{\es}(t_1)v\notin K_\De(M)\right\}\mu^\es(\Dd v)\le\sum_{k=1}^3 \iI_k.
\end{align*}
Combining this with inequalities \ef{0.53}-\ef{0.56}, we arrive at bound \ef{0.55}, where one should replace $\De'$ by $3\De'$.

\subsection{Some a priori bounds}
The following two lemmas are standard, see for instance \cite{temam1979}.
\begin{lemma}\label{0.16}
For any $u_0\in H$ and $\ph\in L^2_{loc}(\rr_+;H)$, we have 
\be\label{0.18}
\|S^\ph(t)u_0-S(t)u_0\|^2\le C e^{c(\|u_0\|^2+t\|h\|^2)}\int_0^t e^{-\lm_1(t-s)}\|\ph(s)\|^2\dd s,
\ee
where inequality holds for all $t\ge0$, and where $C$ and $c$ are positive constants that depend only on the eigenvalue $\lm_1$.
\end{lemma}

\begin{lemma}\label{0.51}
For any positive constant $M$ there is $C(M)>0$ such that for any $T>0$ and $\ph\in L^2(0, T; H)$ satisfying
$$
\int_0^T \|\ph(s)\|^2\dd s\le M,
$$
we have
$$
\sup_{t\in [0, T]}|S^\ph(t)u_0|_{V}\le C(M)\q\text{ for any }u_0\in\aaa.
$$
\end{lemma}

\addcontentsline{toc}{section}{Bibliography}
\bibliography{reff}
\bibliographystyle{plain}

\end{document}